\font\sc=cmcsc10
\font\tencmmib=cmmib10
\font\eightcmmib=cmmib10
\font\tenmsym=msym10
\font\eightmsym=msym8
\font\teneufm=eufm10
\font\tencmsy=cmsy10
\font\eightcmsy=cmsy8
\textfont 7=\tencmsy \scriptfont 7=\eightcmsy 
 
\textfont 9=\tenmsym \scriptfont 9=\eightmsym 
\def\bb {\fam9 } 
\textfont 8=\teneufm \scriptfont 8=\teneufm 
\def\frak {\fam8 } 
\textfont 6=\tencmmib \scriptfont 7=\eightcmmib 
 
\magnification=1200
\font\bgf=cmbx10 scaled\magstep2

\font\tt=cmtt10
 
\centerline{\bgf Massey Products and Deformations}\bigskip

\centerline{\sc Dmitry Fuchs and Lynelle Lang}\smallskip 
 
\centerline{Department of Mathematics, University of California, Davis CA 
95616, USA}\smallskip

\centerline{\tt fuchs@math.ucdavis.edu, lang@math.ucdavis.edu}\bigskip
 
\centerline{\bf 1. Introduction}\medskip
 
It is common knowledge that the construction of one-parameter deformations of 
various algebraic structures, like associative algebras or Lie algebras, 
involves certain conditions on cohomology classes, and that these conditions 
are usually expressed in terms of Massey products, or rather Massey powers. 
The cohomology classes considered are those of certain differential graded Lie 
algebras (DGLAs). It is also known that the Massey products arising in 
deformation theory are not precisely the products considered in the general 
theory of DGLAs (see [R2]). Actually, different natural problems of 
deformation theory give rise to different kinds of Massey products. The 
definitions of these Massey products involve certain equations whose 
coefficients turn out to be, quite unexpectedly, the structure constants of a 
graded commutative associative algebra. Thus to define Massey products in the 
cohomology of a differential graded Lie algebra one should begin by choosing a 
graded commutative associative algebra. It is interesting that, dually, 
different kinds of Massey products in the cohomology of a differential 
commutative associative algebra correspond in a similar way to Lie algebras. 
In particular, the classical Massey products correspond to the Lie algebra of 
strictly upper triangular matrices, while May's matric Massey products [Ma] 
correspond to the Lie algebra of block strictly upper triangular matrices.
 
The article is organized as follows. Section 2 contains a list of various 
Massey-like products which arise in the cohomology of DGLAs. Most of them are 
related to deformations of Lie algebras. In section 3 we give a general 
construction of Massey products in the cohomology of DGLAs. The main result of 
this section is Proposition 3.1, which shows the necessity of the 
associativity of the auxiliary algebra. Section 4 contains an application of 
the construction of Section 3 to Lie algebra deformations, and in Section 5 we 
consider Massey products in the cohomology of graded commutative associative 
differential algebras.
\bigskip

\centerline{\bf 2. Examples.}\medskip
 
{\bf 2.1.} First we recall the standard definition of Massey products in the 
cohomology of a DGLA (see [R2]). A DGLA is a vector space $L$ over a field 
${\bb K}$ of characteristic zero with a ${\bb Z}$ or ${\bb Z}_2$ grading 
$L=\bigoplus_i L^i$ and with a commutator operation
$\mu\colon L\otimes L\to L,\ 
\mu(\alpha,\beta)=[\alpha,\beta]$ of degree 0 and a differential $\delta\colon 
L\to L$ of degree $+1$ satisfying the conditions$$\eqalign{&[\alpha,\beta]=-(-
1)^{\alpha\beta}[\beta,\alpha],\cr &\delta[\alpha,\beta]=[\delta 
\alpha,\beta]+(-1)^\alpha[\alpha,\delta \beta],\cr &[[\alpha,\beta],\gamma]+(-
1)^{\alpha(\beta+\gamma)}[[\beta,\gamma],\alpha]+(-
1)^{\gamma(\alpha+\beta)}[[\gamma,\alpha],\beta]=0,\cr}$$where the degree of a 
homogeneous element is denoted by the same letter as this element. For any 
graded (Lie or co)algebra $A$, define the maps  $S\colon A\otimes A\to 
A\otimes A$ and $C\colon A\otimes A\otimes A\to A\otimes A\otimes A$, by the 
formulas $S(\alpha\otimes \beta)=(-1)^{\alpha\beta}\beta\otimes \alpha,\ 
C(\alpha\otimes\beta\otimes\gamma)=(-
1)^{\alpha(\beta+\gamma)}\beta\otimes\gamma\otimes\alpha$.  Then the 
conditions on $\mu$ and $\delta$ may be rewritten as 
$$\eqalign{&\mu\circ S =-
\mu,\cr &\delta\circ\mu=\mu\circ(\delta\otimes1+1\otimes\delta),\cr 
&\mu\circ(\mu\otimes1)\circ(1+C+C^2)=0.\cr}$$ Note that the tensor product of 
homomorphisms between graded spaces is understood in the graded sense, 
$(f\otimes g)(x\otimes y)=(-1)^{gx}f(x)\otimes g(y)$; in particular, 
$(1\otimes\delta)(\alpha\otimes\beta)=(-1)^\alpha\alpha\otimes\delta(\beta)$. 
 
We denote by $H=\bigoplus_iH^i$ the cohomology of $L$ with respect to the 
differential $\delta$.
 
Let $a_1\in H^{q_1},\dots ,a_r\in H^{q_r}$ be cohomology classes.  If $K$ and
$L$ are proper subsets of $R=\{1,\dots ,r\}$, then define
$$\varepsilon(K,L)=\sum(q_{k}+1)(q_{l}+1),$$ 
 where the sum is taken over 
all $k$ and $l$  such that $k\in K,l\in L$ and $k>l$. Also, for 
any subset $I$ of $R$, let $$P(I) = \left\{ \left(K,L \right) |K,L\subset I, 
K\cup L=I,K\cap L=\emptyset,K\ne\emptyset, L\ne \emptyset\right\}.$$ We say 
that $b\in\langle a_1,\dots ,a_r\rangle$ if there exists a function assigning 
to each proper subset $I=\{i_1,\dots ,i_s\}$ of $R=\{1,\dots ,r\}$ an element 
$\alpha_I\in L^{q_{i_1}+\dots +q_{i_s}-(s-1)}$, such that 
$$\eqalign{&\alpha_i\in a_i,\cr &\delta\alpha_I={1\over2}\sum_{(K,L) \in P(I)} 
(-1)^{\varepsilon(K,L)+\alpha_K+1}[\alpha_K,\alpha_{L}];\cr}$$and$$ 
b\ni{1\over2}\sum_{(K,L) \in P(R)}(-
1)^{\varepsilon(K,L)+\alpha_K+1}[\alpha_K,\alpha_{L}].$$ Hence $\langle 
a_1,\dots  ,a_r\rangle$ is a subset, maybe empty, of $H^{q_1+\dots +q_r-(r-
1)}$. It is  called the Massey product of $a_1,\dots ,a_r$. Note that $\langle 
a_1\rangle =0$ and $\langle a_1,a_2\rangle =[a_1,a_2].$ Also note that 
$\langle a_1,\dots ,a_r\rangle$ is nonempty if and only if for any proper 
subset $\{i_1,\dots ,i_s\}$ of $R$ the Massey product $\langle a_{i_1},\dots 
,a_{i_s}\rangle$ contains zero.\smallskip 
 
{\bf 2.2.} Now we will consider deformations of Lie algebras. Let $\frak g$ be 
a Lie algebra over ${\bb K}$ with bracket [ , ]. A {\it formal 
deformation} of $\frak g$ is defined as a power series 
$$[g,h]_t=[g,h]+\sum_{k=1}^\infty t^k\gamma_k(g,h),\eqno(1)$$which makes 
${\frak g}[[t]]$, the space of formal power series with coefficients in $\frak 
g$, a Lie algebra. This deformed bracket is antisymmetric and satisfies the 
Jacobi identity exactly when$$\gamma_k(g,h)=-
\gamma_k(h,g),$$and$$\delta\gamma_k=-{1\over2}\sum_{i=1}^{k-
1}[\gamma_i,\gamma_{k-i}].\eqno(2)$$Here each $\gamma_k$ is regarded as a 
cochain in
$C^2({\frak g};{\frak g})= {\rm Hom }_k(\bigwedge^2{\frak g}, {\frak g})$,
and the product $[\gamma_i,\gamma_{k-
i}]$ is taken with respect to the standard graded Lie algebra structure in 
$C^\ast({\frak g};{\frak g})$; if $\gamma'\in C^{q'}({\frak g};{\frak 
g}),\gamma''\in C^{q''}({\frak g};{\frak g})$, then $[\gamma',\gamma'']\in 
C^{q'+q''-2}({\frak g};{\frak g})$ (see, for example, [Fu]). If we set 
$C^q({\frak g};{\frak g})=L^{q-1}$, then $L=\bigoplus L^i$ is a DGLA. 
 
For $k=1$ condition (2) means that $\delta\gamma_1=0$ and elements of 
$H^2({\frak g};{\frak g})$ are called {\it infinitesimal deformations} of 
$\frak g$. The cohomology class of $\gamma_1$ is called the {\it differential} 
of the formal deformation (1), or the infinitesimal deformation defined by the 
formal deformation (1). A given infinitesimal deformation $c\in H^2({\frak 
g};{\frak g})$ is the differential of some formal deformation if and only if 
there exist cochains $\gamma_i\in C^2({\frak g};{\frak g})$ such that 
$\delta\gamma_1=0,\gamma_1\in c$, and the whole sequence $\{\gamma_i\}$ 
satisfies the conditions in (2). These conditions are usually referred to as 
$\underbrace{\langle c,\dots ,c\rangle}_r\ni0$ for all $r$; but there exists a 
difference between this $\langle c,\dots ,c\rangle$ and the Massey product 
defined in 2.1. (This difference was pointed out by Retakh in [R1].) Indeed, 
in 2.1 we chose cochains $a_I$ for all $I\subset\{1,\dots ,r\}$, while here we 
take we take just one cochain $\gamma_k$ for each $k$, which means, from the 
point of view of 2.1, that all cochains $a_I$ with sets $I$ of the same 
cardinality are chosen equal (the products are then multiples of each 
other).\smallskip 
 
{\bf 2.3.} Let $\frak g$ be as above. A {\it singular deformation} of $\frak 
g$ is a formal deformation defined as in (1) with $\gamma_1=0$, that is 
$$[g,h]_t=[g,h]+\sum_{k=2}^\infty t^k\gamma_k(g,h).\eqno(3)$$

It is known (see 
[FF]) that some Lie algebras have {\it essentially} singular deformations (those
which cannot be approximated by nonsingular ones). For singular deformations 
the first two conditions from (2) become $$\delta\gamma_2=0,\ 
\delta\gamma_3=0.$$The problem can be stated as follows: for two given 
cohomology classes $c_2,c_3\in H^2({\frak g};{\frak g})$ does there exist a 
singular deformation (3) with $\gamma_2\in c_2,\gamma_3\in c_3$? The equations 
in (2) provide necessary and sufficient conditions for that. The first three 
conditions, $\delta\gamma_4=-{1\over2}[\gamma_2,\gamma_2]$, 
$\delta\gamma_5=-[\gamma_2,\gamma_3]$, and $\delta\gamma_6=-
[\gamma_2,\gamma_4]-{1\over2}[\gamma_3,\gamma_3]$, can be expressed in terms 
of Massey products and powers as $[c_2,c_2]=0$, $[c_2,c_3]=0$, and $- 
{1\over2}[c_3,c_3]\in\langle c_2,c_2,c_2\rangle$.  The fourth condition, 
$\delta\gamma_7=-[\gamma_2,\gamma_5]-[\gamma_3,\gamma_4]$, uses a Massey-like 
product which is simply a restricted form of $\langle c_2,c_2,c_3\rangle$; 
however, the next equation, $\delta\gamma_8= -[\gamma_2,\gamma_6]-
[\gamma_3,\gamma_5]-{1\over2}[\gamma_4,\gamma_4]$, not only is related to a 
combination of two Massey products, $\langle c_2,c_2,c_2,c_2\rangle$ and 
$\langle  c_2,c_3,c_3\rangle$, but these products may both be empty (we do not 
have either the equality $[c_3,c_3]=0$, or the inclusion $\langle  
c_2,c_2,c_2\rangle\ni0$). The subsequent conditions also refer to Massey-like 
products, (also similar to combinations of empty Massey products) which are 
essentially different from those of 2.1 and 2.2.\smallskip 
 
{\bf 2.4.} Let $\frak g$ be as above. Consider two formal deformations, 
$$[g,h]_t=[g,h]+\sum_{k=1}^\infty t^k\gamma_k(g,h),\qquad
[g,h]'_t=[g,h]+\sum_{k=1}^\infty t^k\gamma'_k(g,h)\eqno(4)$$with 
$\gamma'_1=\gamma_1$. For $i\ge2$ set $\beta_i=\gamma'_i-\gamma_i$. The two 
deformations in (4) satisfy the conditions (2), that is$$\delta\gamma_k=-
{1\over2}\sum_{i=1}^{k-1}[\gamma_i,\gamma_{k-i}],\quad \delta\gamma'_k=-
{1\over2}\sum_{i=1}^{k-1}[\gamma'_i,\gamma'_{k-i}].$$If we substitute 
$\gamma_i+\beta_i$ for $\gamma'_i$ and then subtract the first 
equality from the second, we get the following pair of relations: 
$$\delta\gamma_k=-{1\over2}\sum_{i=1}^{k-1}[\gamma_i,\gamma_{k-i}],\quad 
\delta\beta_k=-\sum_{i=1}^{k-2}[\gamma_i,\beta_{k-i}]-{1\over2}\sum_{i=2}^{k-
2}[\beta_i,\beta_{k-i}].\eqno(5)$$In particular, we have $\delta\gamma_1=0$ 
and $\delta\beta_2=0.$ The following problem arises. For two given cohomology 
classes $c,b\in H^2({\frak g};{\frak g})$ does there exist a pair of formal 
deformations (4) with $\gamma'_1=\gamma_1\in c,\ \gamma'_2-\gamma_2=\beta_2\in 
b$? The conditions for this are provided by the relations in (5). First of 
all, the conditions $\delta\gamma_2=-{1\over2}[\gamma_1,\gamma_1]$ and 
$\delta\beta_3=- [\gamma_1,\beta_2]$ mean that $[c,c]=0$ and $[c,b]=0$ in 
order for $\gamma_2$ and $\beta_3$ to exist.  The next pair of relations says 
that $\delta\gamma_3=- [\gamma_1,\gamma_2]$ and $\delta\beta_4 = -
[\gamma_1,\beta_3] - [\gamma_2,\beta_2]-{1\over2}[\beta_2,\beta_2]$, which 
implies that $\langle  c,c,c\rangle\ni0$ and $-{1\over2}[b,b]\in\langle 
c,c,b\rangle$.  The third pair, $\delta\gamma_4 = -[\gamma_1,\gamma_3]-
{1\over2}[\gamma_2,\gamma_2]$ and $\delta\beta_5 = -[\gamma_1,\beta_4]-
[\gamma_2,\beta_3]-[\gamma_3,\beta_2]-[\beta_2,\beta_3]$, may be interpreted 
separately as $\langle c,c,c,c\rangle\ni0$ (in the sense of 2.2) and a product 
related to a combination of $\langle c,c,c,b\rangle$ and $\langle 
c,b,b\rangle$ (empty Massey products) contains zero.  However, the cochains 
$\gamma_2$ and $\gamma_3$ are shared by the equations and so in order for 
$\gamma_4$ and $\beta_5$ to exist, $\langle c,c,c,c\rangle$ must contain zero 
using cochains which also make the second product contain zero.  The 
subsequent relations mean that certain Massey-like products involving $c$ and 
$b$ contain zero with cochains that put $0\in \underbrace{\langle c,\dots 
,c\rangle}_r$ for all $r$ (in the sense of  2.2). These examples show that 
very natural questions about formal deformations of Lie algebras lead to 
different kinds of Massey-like products of cohomology classes in $H^2({\frak 
g};{\frak g})$. In the next section we will show that all these products are 
actually described by a very simple general construction, which involves an 
almost arbitrary auxiliary (graded) associative commutative algebra.
 \bigskip

\centerline{\bf 3. General construction}\medskip
 
{\bf 3.1. Inductive definition.} Let $L=\bigoplus L^q,\ H=\bigoplus H^q$ be as 
above. Suppose that we are given a finite or infinite sequence of integers 
$\{q_k\in{\bb Z}\}_{k=1}^{n\ {\rm or}\ \infty}$ and an indexed set 
$\{c_k^{ij}\in {\bb K}\}_{i,j,k=1}^{n\ {\rm or}\ \infty}$, which satisfy the 
following condition:$$c_k^{ij}=0,\ {\rm if}\ i\ge k, j\ge k\ {\rm or}\  q_k\ne 
q_i+q_j-1.$$ Also we suppose that some integer $r$ is fixed, and $c_k^{ij}=0$ 
for $k\le r$.\smallskip 
 
Let $a_1\in H^{q_1},\dots ,a_r\in H^{q_r}$, and let $b\in H^{q_{r+s}+1}$ for 
some $s\ge1$. We say that $$b\in\langle a_1,\dots ,a_r\rangle_s,$$ if there 
exists a sequence $\alpha_k\in L^{q_k}, k=1,\dots ,r+s-1$, such that\smallskip
 
(i) $\delta\alpha_k=\sum_{i,j}c_k^{ij}[\alpha_i,\alpha_j]$ for $1\le k\le r+s-
1$; in particular, $\delta\alpha_1=0,\dots ,\delta\alpha_r=0$;\smallskip
 
(ii) $\alpha_1\in a_1,\dots ,\alpha_r\in a_r$;\smallskip
 
(iii) $\delta\beta=0$ and $\beta\in b$, where 
$\beta=\sum_{i,j}c_{r+s}^{ij}[\alpha_i,\alpha_j]$.\smallskip 
 
Note that in order for $\langle a_1,\dots ,a_r\rangle_s$ to be  nonempty it is 
necessary that $\langle a_1,\dots ,a_r\rangle_t\ni0$ for all $t<s$.\smallskip
 
All the examples of Section 2 are particular cases of this definition. But in 
all these examples the coefficients $c_k^{ij}$ have two more properties, which 
make the definition particularly convenient.

First, for $i\ne j$, the two terms 
$c_k^{ij}[\alpha_i, \alpha_j]$ and $c_k^{ji}[\alpha_j, \alpha_i]$
in the sum 
$\sum_{i,j}c_k^{ij}[\alpha_i,\alpha_j]$
both contribute the same cochain, 
$[\alpha_i,\alpha_j]=\pm[\alpha_j,\alpha_i]$ 
with possibly different coefficients. It is convenient to make the 
the two terms equal with the 
condition 
$$c_k^{ij}=(-1)^{q_iq_j-1}c_k^{ji}.$$
 
Second, in all the examples above the coefficients $c_k^{ij}$ were chosen in 
such a way that the sum $\sum_{i,j}c_k^{ij}[\alpha_i,\alpha_j]$ is 
automatically a cocycle (belongs to Ker$\, \delta$). 
 
Both of these conditions assume a compact form if we consider the numbers 
$c_k^{ij}$, or rather $(-1)^{q_i-1}c_k^{ij}$, as the structure constants of a 
certain algebra, that is if we consider a ${\bb K}$-algebra $G$ with a basis 
$g^k,k=1,2,\dots$ and with $g^ig^j=\sum_k(-1)^{q_i-1}c_k^{ij}g^k$. Then the 
first

condition on the $c^k_{ij}$ becomes evidently that of graded commutativity (we put 
deg$\, g^k=q_k-1$), while the second condition turns out to be that of 
associativity. To demonstrate this it is convenient to have the previous 
construction in a more general form; in this compact construction we will use 
rather an auxiliary coassociative coalgebra than an associative 
algebra.\smallskip 
 
{\bf 3.2. Direct definition.} Let $L,\mu,\delta,S,C$ and $H$ be the same as in 
2.1. Now let $F$ be a graded cocommutative coassociative coalgebra, that is a 
($\bb Z$ or ${\bb Z}_2$) graded vector space $F$ with a degree 0 mapping 
(comultiplication) $\Delta\colon F\to F\otimes F$ satisfying the conditions 
$S\circ\Delta=\Delta$ and $(1\otimes\Delta)\circ\Delta=(\Delta\otimes1) 
\circ\Delta.$ 
 
Suppose also that in $F$ a filtration $F_0\subset F_1\subset F$ is given, such 
that $F_0\subset{\rm Ker}\, \Delta$ and ${\rm Im}\, \Delta\subset F_1\otimes 
F_1$. The next proposition is, technically, the central fact of this 
article.\smallskip 
 
{\sc Proposition 3.1.} {\it Suppose a linear mapping $\alpha\colon F_1\to L$ 
of degree 1 satisfies the condition$$\delta\circ\alpha=\mu\circ(\alpha\otimes\alpha) 
\circ\Delta. \eqno(6)$$Then$$\mu\circ(\alpha\otimes\alpha)\circ\Delta(F) 
\subset{\rm Ker}\, \delta.$$}\smallskip 
 
{\sc Remark.} The left-hand side of the last formula is well defined because 
$\Delta(F)$ is contained in $F_1\otimes F_1$, the domain of 
$\alpha\otimes\alpha$.\smallskip 
 
Proposition 3.1 will be proved in Subsection 3.3 below.\smallskip
 
{\sc Definition.} Let $a\colon F_0\to H$, $b\colon F/F_1\to H$ be two linear 
maps of degree 1. We say that $b$ is contained in the Massey $F$-product of 
$a$, and write $b\in\langle a\rangle_F$, or $b\in\langle a\rangle$,or 
$b\in\langle a\rangle_{F;F_0,F_1}$, if there exists a degree 1 linear mapping 
$\alpha\colon F_1\to L$ satisfying condition (6), and such that the diagrams 
 
$$ 
\def\mapright#1{\smash{\mathop{\relbar\joinrel\longrightarrow}\limits^{#1}}} 
\def\longmapright#1{\smash{\mathop{\relbar\joinrel\relbar\joinrel\relbar
\joinrel\relbar\joinrel\relbar\joinrel\rightarrow}\limits^{#1}}} 
\def\mapdown#1{\Big\downarrow\rlap{$\vcenter{\hbox{$\scriptstyle#1$}}$}} 
\matrix{F_0&\mapright{\alpha|_{F_0}}&{\rm Ker}\, \delta\cr 
\Big\Vert&&\mapdown{\pi}\cr F_0&\mapright{a}&H\cr},\qquad\qquad 
\matrix{F&\longmapright{\mu\circ(\alpha\otimes\alpha)\circ\Delta}&{\rm Ker}\, 
\delta\cr \mapdown{\pi}&&\mapdown{\pi}\cr F/F_1&\longmapright{b}&H\cr}
\eqno(7)$$are commutative, where the vertical maps labeled $\pi$ denote the 
projections of each space onto the quotient space. 
 
Note that the upper horizontal maps of the diagrams in (7) are well defined, 
since $\alpha(F_0)\subset\alpha({\rm Ker}\, \Delta)\subset{\rm Ker}\, \delta$ 
by virtue of (6), and $\mu\circ(\alpha\otimes\alpha)\circ\Delta(F) \subset{\rm 
Ker}\, \delta$ by Proposition 3.1. 
 
The inductive definition of 3.1 is a particular case of this direct definition 
in the following way. One can take for $F$ the vector space spanned by 
$f_1,\dots,f_{r+s}$, define $F_0$ and $F_1$ as subspaces spanned by $f_j$ with 
$j\le r$ and $j<r+s$ respectively, and let deg$\, f_k=q_k-1$, $\Delta 
f_k=\sum_{i,j}(-1)^{q_i-1}c^{ij}_kf_i\otimes f_j$. The maps $a\colon F_0\to H$ 
and $b\colon F/F_1\to H$ are defined by the formulas $f_k\mapsto a_k\ (k\le 
r)$ and $f_{r+s}+F_1\mapsto b$. The map $\alpha$ takes $f_k$ into $\alpha_k$; 
condition (6) and the commutativity of diagrams (7) become conditions (i),
(ii), and (iii) of 3.1.\smallskip 
 
Note in conclusion that our definition makes sense even in the case, when 
$F_1=F$. In this case we do not need to specify any $b$, and we will simply 
say that $a$ {\it satisfies the condition of triviality of Massey 
$F$-products}. \smallskip 
 
{\bf 3.3. Proof of Proposition 3.1.} We need to prove 
that$$\delta\circ\mu\circ (\alpha\otimes\alpha)\circ\Delta=0.$$The definition 
of a DGLA and condition (6) 
imply:$$\eqalign{\delta\circ\mu\circ(\alpha\otimes\alpha)\circ\Delta&= 
\mu\circ (\delta\otimes1+1\otimes\delta)\circ(\alpha\otimes 
\alpha)\circ\Delta\cr &=\mu\circ((\delta\circ\alpha)\otimes\alpha-
\alpha\otimes(\delta\circ\alpha))\circ\Delta\cr 
&=\mu\circ((\delta\circ\alpha)\otimes\alpha)\circ\Delta-
\mu\circ(\alpha\otimes(\delta\circ\alpha))\circ\Delta;\cr 
\mu\circ((\delta\circ\alpha)\otimes\alpha)\circ\Delta 
&=\mu\circ((\mu\circ(\alpha\otimes\alpha)\circ\Delta) 
\otimes\alpha)\circ\Delta\cr 
&=\mu\circ((\mu\circ(\alpha\otimes\alpha)\circ\Delta) 
\otimes(1\circ\alpha\circ1))\circ\Delta\cr 
&=\mu\circ(\mu\otimes1)\circ(\alpha\otimes\alpha\otimes\alpha) 
\circ(\Delta\otimes1)\circ\Delta;\cr}$$and 
similarly$$\mu\circ(\alpha\otimes(\delta\circ\alpha))\circ\Delta= 
\mu\circ(1\otimes\mu)\circ(\alpha\otimes\alpha\otimes\alpha)\circ 
(1\otimes\Delta)\circ\Delta.$$ 
 
The following two lemmas hold even if the coalgebra $F$ is neither 
cocommutative, nor coassociative.\smallskip 
 
{\sc Lemma 1.} $$\Delta\otimes1 = C\circ(1\otimes\Delta)\circ  S\colon 
F\otimes F\to L\otimes L\otimes L.$$ 
 
{\sc Proof}.  Let $f,g \in F.$  Then $$\eqalign {C\circ(1\otimes\Delta )\circ 
S(f\otimes g) &= (-1)^{fg}C\circ (1\otimes\Delta )(g\otimes f) \cr &= (-
1)^{fg}C\circ g\otimes\Delta(f) \cr&= \Delta(f)\otimes g \cr&= (\Delta\otimes 
1)(f\otimes g).}$$\smallskip 
 
{\sc Lemma 2.} $$C\circ(\alpha\otimes\alpha\otimes\alpha) = 
(\alpha\otimes\alpha\otimes\alpha)\circ C\colon F\otimes F\otimes F \to 
L\otimes L\otimes L. $$ 
 
{\sc Proof}.  Let $f,g,h \in F.$ Then  $$\eqalign{ 
C\circ(\alpha\otimes\alpha\otimes\alpha)(f\otimes g\otimes h)& = (-1)^gC\circ 
\alpha(f)\otimes \alpha(g)\otimes \alpha(h) \cr &= (-1)^{g + 
(f+1)(g+h)}\alpha(g)\otimes \alpha(h)\otimes \alpha(f),\cr 
(\alpha\otimes\alpha\otimes\alpha)\circ C(f\otimes g\otimes h)&= (-
1)^{f(g+h)}(\alpha\otimes\alpha\otimes\alpha)(g\otimes h\otimes f)\cr &= (-
1)^{f(g+h) +h}\alpha(g)\otimes \alpha(h)\otimes \alpha(f).\cr}$$ 
 
Using the Lemmas, the cocommutativity of $F$, and the coassociativity of 
$F$:$$\eqalign{\mu\circ(\mu\otimes1)\circ(\alpha\otimes\alpha\otimes\alpha)
&\circ(\Delta\otimes1)\circ\Delta\cr 
&=\mu\circ(\mu\otimes1)\circ(\alpha\otimes\alpha\otimes\alpha)\circ 
C\circ(1\circ\Delta)\circ 
S\circ\Delta\cr &=\mu\circ(\mu\otimes1)\circ 
C\circ(\alpha\otimes\alpha\otimes\alpha)\circ(1\circ\Delta)\circ 
S\circ\Delta\cr &=\mu\circ(\mu\otimes1)\circ 
C\circ(\alpha\otimes\alpha\otimes\alpha)\circ(1\circ\Delta)\circ\Delta\cr 
&=\mu\circ(\mu\otimes1)\circ 
C\circ(\alpha\otimes\alpha\otimes\alpha)\circ(\Delta\otimes1)\circ\Delta.
\cr}$$ Similarly,$$\eqalign{\mu\circ(\mu\otimes1)\circ 
C\circ(\alpha\otimes\alpha\otimes\alpha)&\circ(\Delta\otimes1)\circ\Delta\cr 
&=\mu\circ(\mu\otimes1)\circ 
C^2\circ(\alpha\otimes\alpha\otimes\alpha)\circ(\Delta\otimes1)\circ\Delta.
\cr}$$Hence, 
$$\eqalign{\mu\circ(\mu\otimes1)\circ(\alpha\otimes\alpha\otimes\alpha)
&\circ(\Delta\otimes1)\circ\Delta\cr &={1\over3} 
\mu\circ(\mu\otimes1)\circ(1+C+C^2)\circ(\alpha\otimes\alpha\otimes\alpha) 
\circ(\Delta\otimes1)\circ\Delta=0.\cr}$$In the same way,
$$\mu\circ(1\otimes\mu)\circ(\alpha\otimes\alpha\otimes\alpha)\circ(1\otimes 
\Delta)\circ\Delta=0,$$and hence, 
$$\delta\circ\mu\circ(\alpha\otimes\alpha)\circ\Delta=0.$$\smallskip
 
{\bf 3.4. Survey of examples of Section 2.} Our goal is to show that all the 
Massey-like products considered in Section 2 are covered by the previous 
construction. 
 
For the definition of Massey products in 2.1 one should set$$F={\rm 
span}\{f_I|\emptyset\ne I\subseteq R=\{1,\dots ,r\}\}$$$$F_0={\rm 
span}\{f_{\{1\}},\dots ,f_{\{r\}}\},\ F_1={\rm span}\{f_I|I\ne R\},$$$${\rm 
deg}\, f_I=\sum_{i\in I}(q_i-1),$$$$\Delta f_I={1\over2}\sum_{K\subset I\atop 
\emptyset\ne K\ne I}(-1)^{\varepsilon(K,I-K)}f_K\otimes f_{I-K}.$$ The dual 
algebra $G=F^\ast$ is spanned by $g^I,\emptyset\ne I\subseteq R$, and the 
multiplication in $G$ is defined by the formula$$g^Kg^L=\cases{0&if $K\cap 
L\ne\emptyset$,\cr {1\over 2}(-1)^{\varepsilon(K,L)}g^{K\cup L}&if $K\cap 
L=\emptyset$.\cr}$$If we put $x_i=g^{\{i\}}$, then for $I=\{i_1,\dots 
,i_s\}$$$g^I=2^{1-s}x_{i_1}\dots x_{i_s}.$$Hence $G$ is the maximal ideal in 
the quotient of the free graded commutative algebra with the generators 
$x_1,\dots ,x_r$ with deg$\, x_i=q_i-1$ over the ideal generated by 
$x_1^2,\dots ,x_r^2.$ By the way, $F_0^\ast=G/G^2$, and $F_1^\ast=G/G^r$. 
 
For the condition in 2.2 that $\underbrace{\langle c,\dots ,c\rangle}_r\ni0$ 
for all $r$ one should set$$F={\rm span}\{f_1,f_2,\dots\},$$$$F_0={\bb K}f_1,\ 
F_1=F,$$$${\rm deg}\, f_i =0,$$$$\Delta f_i=-{1\over2}\sum_{k=1}^{i-
1}f_k\otimes f_{i-k}.$$The dual algebra $G$ is formed by finite or infinite 
linear combinations of $g^k, k=1,2,\dots$ with the multiplication$$g^kg^l=-
{1\over2}g^{k+l}.$$If we put $g^1=t$, then $g^k=(-2)^{1-k}t^k$ (in the left 
hand side $k$ is a superscript, while in the right hand side it is the 
exponent). Hence $G$ is the maximal ideal in ${\bb K}[[t]]$, the algebra of 
formal power series in one variable $t$ of degree 0. Again $F_0^\ast=G/G^2$. 
 
For the definition of 2.3 one has $$F={\rm span}\{f_2,f_3,\dots\},$$and all 
the formulas are obvious modifications of the formulas of the previous 
paragraph. In particular, $G$ is the square of the maximal ideal of ${\bb 
K}[[t]]$. If we set $u=t^2, v=t^3$, we may interpret $G$ as the maximal ideal 
in ${\bb K}[[u,v]]/(u^3-v^2)$. Again $F_0^\ast=G/G^2$. 
 
For the problem in 2.4 of the existence of a pair of deformations one should 
take$$F={\rm span}\{f_1,f_2,f_3,\dots;\varphi_2,\varphi_3,\dots\},$$$$F_0={\rm 
span}\{f_1,\varphi_2\},\ F_1=F,$$$${\rm deg}(f_i)=0,{\rm 
deg}(\varphi_i)=0,$$$$ \Delta f_i=-{1\over2}\sum_{k=1}^{i-1}f_k\otimes f_{i-
k},\ \Delta\varphi_i=-{1\over2}\sum_{k=1}^{i-2}f_k\otimes\varphi_{i-k}-
{1\over2}\sum_{k=2}^{i-2}\varphi_k\otimes\varphi_{i-k}.$$The dual algebra $G$ 
is formed by finite or infinite linear combinations of $g^k, k=1,2,\dots$, and 
$\psi^k, k=2,3,\dots$ with the multiplication$$g^kg^l=-{1\over2}g^{k+l},\ 
g^k\psi^l=-{1\over2}\psi^{k+l},\ \psi^k\psi^l=-{1\over2}\psi^{k+l}.$$If we put 
$g^1=t,\psi^2=u,$ then $g^k=(-2)^{1-k}t^k,\psi^l=(-2)^{2-l}t^{l-2}u$, and $-
2u^2=t^2u$. Hence $G$ is the maximal ideal in ${\bb K}[[t,u]]/(2u^2+t^2u)$, 
and again we have $F_0^\ast=G/G^2$.\bigskip 
 
\centerline{\bf 4. An application to Lie algebra deformations}\medskip
 
{\bf 4.1. Deformations.} We consider the algebraic version of deformation 
theory (see, for example, [Fi]). Let $S$ be a commutative associative $\bb 
K$-algebra with an identity element, and with a distinguished maximal ideal 
${\frak m}\subset S$ with $S/{\frak m}\cong\bb K$; let $\varepsilon\colon S\to 
S/{\frak m}=\bb K$ be the projection with $\varepsilon(1)=1$. Consider a Lie 
algebra $\frak g$ over 
$S$. A {\it deformation} of $\frak g$ with  base $S$ 
is, by definition, a structure of a Lie algebra over 

$\bb K$ on the $S$-module 
${\frak g}\otimes S$, such that $1\otimes\varepsilon\colon{\frak g}\otimes 
S\to{\frak g}\otimes{\bb K}={\frak g}$ is a Lie algebra 
homomorphism.\smallskip 
 
{\sc Example.} Let dim$\, {\frak g}<\infty, {\bb K}={\bb R}, S=C^\infty X,$ 
the algebra of smooth functions on a smooth compact manifold $X$, and $\frak 
m$ be the ideal of functions vanishing at a point $x_0\in X$. Then a 
deformation of $\frak g$ with base $S$ is the same as a smooth family of 
Lie algebra structures on $\frak g$, parametrized by  points $x\in X$ and 
reducing to the initial Lie algebra structure on $\frak g$ for 
$x=x_0$.\smallskip 
 
Suppose that dim$\, S<\infty$. Let $\tau\colon({\frak g}\otimes S)\times({\frak 
g}\otimes S)\to{\frak g}\otimes S$ be a skew-symmetric $S$-bilinear binary 
operation on ${\frak g}\otimes S$ which satisfies the condition of 
$1\otimes\varepsilon$ being a homomorphism, that is 
$(1\otimes\varepsilon)\tau(g\otimes s,h\otimes t)=[g,h]\varepsilon(st)$, but 
not necessarily satisfying the Jacobi identity. For a linear functional 
$\varphi\colon{\frak m}\to{\bb K}$ define a map $\alpha_\varphi\colon{\frak 
g}\otimes{\frak g}\to\frak g$ by the 
formula
$$\alpha_\varphi(g,h)=(1\otimes\varphi)(\tau(g\otimes1,h\otimes1)-
[g,h]\otimes1)$$(it follows from the fact that $1\otimes\varepsilon$
is a homomorphism that 
$\tau(g\otimes1,h\otimes1)-[g,h]\otimes1\in{\frak g}\otimes{\frak m}$). 
Evidently, $\alpha_\varphi$ is bilinear and skew-symmetric, which means that 
$\alpha_\varphi\in C^2({\frak g}\, ;\, {\frak g})$. Moreover, if we
put ${\frak m}^\ast=F$, then $\varphi\mapsto\alpha_\varphi$ is a linear map $\alpha\colon 
F\to C^2({\frak g}\, ;\, {\frak g})$, and it is clear that $\tau$ with the 
properties listed above and $\alpha$ determine each other.
 
Let $\Delta\colon F\to F\otimes F$ be the comultiplication in $F$ dual to the
multiplication in $\frak m$, $\delta\colon C^2({\frak g}\, ;\, {\frak g})\to
C^3({\frak g}\, ;\, {\frak g})$ be the Lie algebra cochain differential,
$\mu\colon C^2({\frak g}\, ;\, {\frak g})\otimes C^2({\frak g}\, ;\,
{\frak g}) \to C^2({\frak g}\, ;\, {\frak g})$ be the Lie product.\smallskip
 
{\sc Proposition 4.1}. {\it The operator $\tau$ satisfies the Jacobi identity 
if and only if $\alpha$ satisfies the Maurer-Cartan equation} 
$$\delta\circ\alpha+{1\over2}\mu\circ(\alpha\otimes\alpha)\circ\Delta=0.$$
\smallskip 
 
This proposition may be regarded as well known, but we will prove it below 
(see 4.3) for the sake of completeness.
 
Let $F_0=({\frak m}/{\frak m}^2)^\ast\subset F$. This is the {\it tangent
space of $S$ at $\frak m$}.  Obviously, $F_0\subset{\rm Ker }\quad \Delta,$
and Proposition 4.1 shows that if $\tau$ is a deformation of $\frak g$, that
is if $\tau$ satisfies the Jacobi identity, then $\alpha(F_0)\subset{\rm
Ker}\, \delta.$ Consider the composition$$a\colon F_0\,
{\buildrel\alpha\over\longrightarrow}\, {\rm Ker}\, \delta\,
{\buildrel{\pi}\over\longrightarrow}\, H^2({\frak g}\, ;\, {\frak
g}),\eqno(8)$$
where $\pi$ is the projection map.An arbitrary linear map $F_0\to H^2({\frak g}\, ;\, {\frak g})$
is called an {\it infinitesimal deformation of $\frak g$ with base $S$}, and
$a$ is called the {\it infinitesimal deformation} determined by the deformation
$\tau$, or the {\it differential} of the deformation $\tau$. Proposition 4.1
implies\smallskip
 
{\sc Theorem 4.2.} {\it An infinitesimal deformation $a\colon F_0\to 
H^2({\frak g}\, ;\, {\frak g})$ is a differential of some deformation with 
base $S$ if and only if $-{1\over2}a$ satisfies the condition of triviality of 
Massey $S$-products}. \smallskip 
 
{\bf 4.2. Formal deformations.} Let $S$ be a local commutative associative 
algebra over $\bb K$ with maximal ideal $\frak m$ and canonical 
projection $\varepsilon\colon S\to\bb K$. Suppose that dim$\, {\frak m}^{k-
1}/{\frak m}^k<\infty$ for all $k$. For a vector space $V$ we define $V\, 
\widehat\otimes\,  S$ as 
${{\rm lim} \atop \longleftarrow} (V\otimes(S/{\frak m}^k))$.
 A {\it formal deformation} of a Lie algebra $\frak g$ with base $S$ is 
defined in the same way as a deformation of $\frak g$ with  base $S$, but 
with ${\frak g}\, \widehat\otimes\,  S$ instead of ${\frak g}\otimes 
S$.\smallskip 
 
{\sc Example.} Let $S={\bb K}[[t]]$. Then a formal deformation of $\frak g$ 
with base $S$ is the same as a 1-parameter formal deformation of $\frak g$ 
(see 2.2).\smallskip
 
Everything that was said in 4.1, including Proposition 4.1 and Theorem 4.2 may be 
repeated for formal deformations with two changes: a functional 
$\varphi\colon{\frak m}\to\bb K$ is assumed to be continuous, that is 
$\varphi({\frak m}^k)=0$ for some $k$; the space $F$ is the space of all {\it 
continuous} functionals ${\frak m}\to\bb K$.
 
Note that the formal version of Theorem 4.2 generalizes the results of 2.2, 
2.3, and 2.4. In these cases $S$ corresponds to the straight line $({\bb 
K}[[t]])$, the semicubic parabola $({\bb K}[[u,v]]/(u^3-v^2))$ and the union 
of two smooth curves with a first order tangency $({\bb 
K}[[t,u]]/(u(t^2+2u)))$.\smallskip 
 
{\bf 4.3. Proof of Proposition 4.1.} We prove here Proposition 4.1 as it was 
stated in 4.1. The formal version from 4.2 is proved in a similar way.
 
Let $\{m_i\}$ be a basis of $\frak m$. Suppose that 
$\tau\colon({\frak 
g}\otimes S)\times({\frak g}\otimes S)\to{\frak g}\otimes S$ satisfies the 
conditions listed in 4.1. Then for $g,h\in\frak 
g$
$$\tau(g\otimes1,h\otimes1)=[g,h]\otimes1+ \sum_i\alpha_i(g,h)\otimes 
m_i$$with some 
 
$\alpha_i\in C^2({\frak g}\, ;\, {\frak g})$. Hence for 
$g,h,k\in\frak g$$$\eqalign{\tau(\tau(g\otimes1,&h\otimes1),k\otimes1) 
=\tau([g,h]\otimes1+ \sum_i\alpha_i(g,h)\otimes m_i,k\otimes1)\cr 
&=[[g,h],k]\otimes1+\sum_i \alpha_i([g,h],k)\otimes 
m_i+\sum_i[\alpha_i(g,h),k]\otimes m_i\cr &\hskip1in 
+\sum_{i,j}\alpha_j(\alpha_i(g,h),k)\otimes m_im_j,\cr}\eqno(10)$$ and the 
Jacobi identity for $\tau$ means that the cyclic sum of this expression with 
respect to $g,h,k$ is equal to 0. 
 
Let $\varphi\colon{\frak m}\to\bb K$ be a linear functional, and let 
$\varphi(m_i)=\varphi_i\in\bb K$. Then $\alpha_\varphi=\sum\varphi_i\alpha_i$.
 
Let $\Delta\varphi=\sum_p\psi'_p\otimes\psi''_p$; we suppose that this sum is 
symmetric, that is with each $\psi'_p\otimes\psi''_p$ it contains also 
$\psi''_p\otimes\psi'_p$. Let $\psi'_p(m_i)=\psi'_{p,i},\, 
\psi''_p(m_i)=\psi''_{p,i}$. Then$$\varphi(m_im_j)= \Delta\varphi(m_i\otimes 
m_j)=\sum_p\psi'_p(m_i)\psi''_p(m_j)=\sum_p \psi'_{p,i}\psi''_{p,j}$$ 
 
Apply $\varphi$ to the last term of 
(10):$$\eqalign{\varphi\left(\sum_{i,j}\alpha_j (\alpha_i(g,h),k)\otimes 
m_im_j\right)&=\sum_{i,j}\varphi(m_im_j)\alpha_j 
(\alpha_i(g,h),k)\cr &=\sum_{i,j,p}\psi'_{p,i}\psi''_{p,j}\alpha_j 
(\alpha_i(g,h),k)\cr &=\sum_p\alpha_{\psi'_p}(\alpha_{\psi''_p}(g,h),k).}$$If 
we apply $\varphi$ to the whole right hand side of (10), we 
get$$\alpha_\varphi([g,h],k)+[\alpha_\varphi(g,h),k]+\sum_p\alpha_{\psi'_p} 
(\alpha_{\psi''_p}(g,h),k).$$The cyclic sum of the last expression is equal to 
$$\delta\alpha_\varphi(g,h,k)+{1\over2}\sum_p[\alpha_{\psi'_p},
\alpha_{\psi''_p}](g,h,k),$$ and the Jacobi identity for $\tau$ means 
precisely that this is equal to 0. Proposition 4.1 follows.\bigskip
 
\centerline{\bf 5. Classical Massey products.}\medskip
 
{\bf 5.1. Definitions.} Originally Massey products were considered not for the 
cohomology of DGLAs, but rather for the cohomology of graded commutative
multiplicative complexes. Let $A=\bigoplus A^i$ be a graded differential 
associative commutative algebra with the differential $\delta$ of degree $+1$, 
and let $H=\bigoplus H^i$ be its cohomology. Let $a_1\in H^{q_1},\dots ,a_r\in 
H^{q_r},b\in H^{q_1+\dots+q_r-(r-1)}$. We say that $b$ belongs to the Massey 
product of $a_1,\dots,a_r$, and write $b\in\langle a_1,\dots,a_r\rangle$ if 
there exist a family $\alpha_{ij}\in A^{q_i+\dots+q_{j-1}-(j-i-1)},\, 1\le 
i<j\le r+1,\, (i,j)\ne(1,r+1)$, such that\smallskip
 
(i) $\delta\alpha_{ij}=\displaystyle{\sum_{k=i+1}^{j-1}(-
1)^{\alpha_{ik}}\alpha_{ik}\alpha_{kj}}$; in particular, 
$\delta\alpha_{i,i+1}=0$; 
 
(ii) $\alpha_{i,i+1}\in a_i$;
 
(iii) $b\ni\displaystyle{\sum_{k=2}^{r}(-
1)^{\alpha_{1k}}\alpha_{1k}\alpha_{k,r+1}}$. \smallskip 
 
What is important here, is that the right hand side of the equality (i) 
belongs to Ker$\, \delta$ if (i) holds for all $\delta\alpha_{i'j'}$ with $j'-
i'<j-i$; in particular, (i) implies that the right hand side of (iii) belongs 
to Ker$\, \delta$.
 
Note that $\langle a\rangle=0,\, \langle a,b\rangle=ab$, and $\langle 
a_1,\dots,a_r\rangle$ is non-empty if and only if $\langle 
a_i,a_{i+1},$ $\dots,a_j\rangle\ni0$ for all $i,j$ such that $1\le i<j\le r,\, 
(i,j)\ne (1,r)$.\smallskip 
 
{\bf 5.2. May's matric Massey products.} The following generalization of the 
previous construction belongs to May [Ma]. Let $p_1,\dots,p_{r+1}$ be 
positive integers, and let $a_i$ denote not an element of $H^{q_i}$, but 
rather a $p_i\times p_{i+1}$ matrix whose entries are elements of $H^{q_i}$. 
Then the definition of 5.1 may be repeated with $b$ being a $p_1\times 
p_{r+1}$ matrix with entries in $H^{q_1+\dots+q_r-(r-1)}$, $\alpha_{ij}$ being 
a $p_i\times p_j$ matrix with entries in $A^{q_i+\dots+q_{j-1}-(j-i-1)}$, and 
the products in (i) and (iii) being considered in the matrix sense. The 
resulting product is called the {\it matric Massey product}. \smallskip 
 
{\bf 5.3. The general construction.} Let $A$ and $H$ denote the same as 
before. Let $Q$ be a graded Lie coalgebra over $\bb K$ with the 
comultiplication $\Delta\colon Q\to Q\otimes Q$, and let $Q_0\subset 
Q_1\subset Q$ be a filtration with $\Delta(Q_0)=0,\, \Delta(Q)\subset 
Q_1\otimes Q_1$. Let $a\colon Q_0\in H,\, b\colon Q/Q_1\to H$ be linear maps 
of degree 1 such that $\delta\circ\alpha=\mu\circ(\alpha\otimes\alpha) 
\circ\Delta$ (where $\mu\colon A\otimes A\to A$ and $\delta\colon A\to A$) are 
the multiplication and the differential in $A$), such that the diagrams 
similar to (7) (with $Q$ and $A$ instead of $F$ and $L$) are commutative.
 
It is important for this definition that the fact similar to Proposition 3.1 
holds: $Q$ being a Lie coalgebra implies that 
$\delta\circ\mu\circ(\alpha\otimes\alpha)\circ\Delta =0$. The proof of this 
fact is a replica of that of Proposition 3.1.
 
For the classical Massey product the corresponding Lie coalgebra is spanned by 
$f_{ij}\ (1\le i\le j\le r+1),\ {\rm deg}\, f_{ij}={\rm deg}\, \alpha_{ij}-1$, 
and $\Delta$ is defined by the formula $\Delta 
f_{ij}=\displaystyle{{1\over2}\sum_{1<k<j}}(f_{ik}f_{kj}-(-
1)^{f_{ik}f_{kj}}f_{kj}f_{ik})$. If we ignore the coefficient $1\over2$ (which 
is the matter of the substitution $f_{ij}=2f'_{ij}$), then the dual graded 
Lie algebra becomes the graded Lie algebra of strictly upper triangular 
$(r+1)\times(r+1)$ matrices.
 
Similarly, the matric Massey product of 5.2 corresponds to the graded Lie 
algebra of block strictly upper triangular matrices with the block sizes 
$p_1,\dots,p_{r+1}$.\bigskip

The authors would like to thank Michael Penkava for discussing these
results with them and providing valuable suggestions.
\bigskip

\noindent {\bf References}\medskip
 
\item{[FF]} A. Fialowski and D. Fuchs, Deformations of the infinite 
dimensional Lie algebra of vector fields $L_1$, {\it in preparation}.
 
\item{[Fi]} A. Fialowski, Deformations of Lie algebras, Math USSR--Sbornik, 
55:2 (1986), 467--473.
 
\item{[Fu]} D. Fuchs, Cohomology of infinite dimensional Lie algebras 
(Consultants Bureau, New York, 1986).
 
\item{[Ma]} J. P. May, Matric Massey products, J. of Algebra 12 (1969) 
533--568.
 
\item{[R1]} V. Retakh, The Massey operations in Lie superalgebras and 
deformations of complex-analytic algebras, Funct. Anal. and Appl. 11 (1977) 
319--321. 
 
\item {[R2]} V. Retakh, Lie-Massey brackets and $n$-homotopically 
multiplicative maps of differential graded Lie algebras, J. Pure and Appl. 
Algebra 89 (1993) 217--229.\bye